\documentclass[reqno,tbtags]{amsproc}
\usepackage{amssymb}
\usepackage{epsfig}
\usepackage{euscript}
\usepackage{eufrak}
\newtheorem{theorem}{Theorem}[section]
\newtheorem{corollary}{Corollary}[section]
\theoremstyle{definition}
\newtheorem{remark}{Remark}[section]
\def\AAS{{Annals of Actuarial Science\/}}
\def\SAT{{Skan\-di\-na\-visk Ak\-tua\-rie\-tid\-skrift\/}}
\def\IME{{In\-su\-ran\-ce: Ma\-the\-ma\-tics and Eco\-no\-mics\/}}
\def\SAJ{{Scan\-di\-na\-vian Ac\-tua\-rial Journal\/}}
\def\AB{{ASTIN Bulletin\/}}
\def\TPA{{Theory Probab. Appl.}}
\def\AAP{{Ad\-van\-ces in Ap\-pli\-ed Pro\-ba\-bi\-li\-ty\/}}
\def\JAP{{Journal of Applied Probability\/}}
\def\SPA{{Stochastic Processes and their Applications\/}}

\numberwithin{equation}{section}
\renewcommand{\P}{\mathsf{P}}
\newcommand{\D}{\mathsf{D}}
\newcommand{\E}{\mathsf{E}}
\newcommand{\p}{\mathsf{p}}
\newcommand{\AInt}[3]{{\mathcal{#1}}_{#2}(#3)}
\newcommand{\X}[1]{X_{#1}}
\newcommand{\Sum}[1]{S_{#1}}
\newcommand{\Ass}[1]{\bar{#1}}
\newcommand{\adjustL}{\varkappa}
\newcommand{\moplus}{m_{_{\scriptscriptstyle\vartriangle}}}
\newcommand{\mominus}{m_{_{\scriptscriptstyle\triangledown}}}
\newcommand{\Doplus}{D_{{\scriptscriptstyle\vartriangle}}}
\newcommand{\Dominus}{D_{{\scriptscriptstyle\triangledown}}}
\newcommand{\Do}{Q}
\newcommand{\hv}[1]{\mathsf{v}_{#1}}
\newcommand{\RPtimeR}[1]{\Upsilon^{\text{\rm [ren]}}_{#1\,\mid\,\paramT,\paramY}}
\newcommand{\BesselI}[1]{I_{#1}}
\def\paramY{\varrho}
\def\paramT{\delta}
\newcommand{\RtimeR}[1]{\Upsilon^{\text{\rm [ren]}}_{#1}}
\newcommand{\Y}[1]{Y_{#1}}
\newcommand{\T}[1]{T_{#1}}
\newcommand{\homN}[1]{N_{#1}}
\newcommand{\homV}[1]{V_{#1}}
\newcommand{\difP}{\vartheta}
\newcommand{\W}[1]{\mathsf{W}_{#1}}
\newcommand{\DPtimeR}[1]{\Upsilon^{\text{\rm [dif]}}_{#1\,\mid\,\difP}}
\newcommand{\UGauss}[2]{\varPhi_{\left({#1},{#2}\right)}}
\newcommand{\Ugauss}[2]{\varphi_{\left({#1},{#2}\right)}}
\newcommand{\muIG}{\mu}
\newcommand{\HmuIG}{\hat{\mu}}
\newcommand{\lambdaIG}{\lambda}
\newcommand{\cS}{c^{*}}

\begin{document}
%%%%%%%%%%%%%%%%%%%%%%%%%%%%%%%%%%%%%%%%%%%%%%%%%%%%%%%%%%%%%%%%%%%%%%%%%%%%%%%%%%%%%%%%
\author[Vsevolod K. Malinovskii]{\Large Vsevolod K. Malinovskii}


\address{Central Economics and Mathematics Institute (CEMI) of Russian Academy of Science,
117418, Nakhimovskiy prosp., 47, Moscow, Russia}

\email{Vsevolod.Malinovskii@mail.ru, admin@actlab.ru}

\title[ON THE TIME OF FIRST LEVEL CROSSING]{ON APPROXIMATIONS FOR THE DISTRIBUTION OF
THE TIME OF FIRST LEVEL CROSSING}

\maketitle

%%%%%%%%%%%%%%%%%%%%%%%%%%%%%%%%%%%%%%%%%%%%%%%%%%%%%%%%%%%%%%%%%%%%%%%%%%%%%%%%%%%%%%%%
\section{Introduction}\label{ewrtrjkr}
%%%%%%%%%%%%%%%%%%%%%%%%%%%%%%%%%%%%%%%%%%%%%%%%%%%%%%%%%%%%%%%%%%%%%%%%%%%%%%%%%%%%%%%%

This paper is an overview of the classical level crossing problem which is
studied extensively in the literature and is fundamental in many branches of
applied probability. We discuss a number of approximations with an emphasis on
their performance, methods of justification and technical conditions which are
required in these methods, including a new approximation called ``inverse
Gaussian''. It is derived by a new method, is fruitful for solving related
problems, and is valid under mild regularity conditions. We emphasize its
novelty and boons.

%%%%%%%%%%%%%%%%%%%%%%%%%%%%%%%%%%%%%%%%%%%%%%%%%%%%%%%%%%%%%%%%%%%%%%%%%%%%%%%%%%%%%%%%
\subsection{Level crossing by a diffusion process}\label{wrethjrtkj}
%%%%%%%%%%%%%%%%%%%%%%%%%%%%%%%%%%%%%%%%%%%%%%%%%%%%%%%%%%%%%%%%%%%%%%%%%%%%%%%%%%%%%%%%

Let $\homV{s}=\difP s+\sigma\W{s}$, $s>0$, be a diffusion process with drift
coefficient $\difP>0$ and diffusion coefficient $\sigma>0$, and let $c$ be a
positive constant. The random variable
\begin{equation*}
\DPtimeR{u,c}=\inf\big\{s>0:(\difP-c)s+\sigma\W{s}>u\big\},
\end{equation*}
or $+\infty$, if $(\difP-c)s+\sigma\W{s}\leqslant u$ for all $s>0$, is the
first passage time to level $u>0$ of the shifted diffusion process
$\homV{s}-cs$, $s>0$. It is well known that with $\UGauss{0}{1}(x)$ denoting
compound distribution function (c.d.f.) of a standard Gaussian distribution, we
have
\begin{multline*}
\P\big\{\DPtimeR{u,c}\leqslant t\big\}=\P\Big\{\sup_{0\leqslant s\leqslant
t}\big((\difP-c)s+\sigma\W{s}\big)>u\Big\}
\\[0pt]
=1-\UGauss{0}{1}\left(\frac{u-(\difP-c)t}{\sigma\sqrt{t}}\right)
+\exp\left\{2\frac{(\difP-c)
u}{\sigma^2}\right\}\,\UGauss{0}{1}\left(\frac{-u-(\difP-c)
t}{\sigma\sqrt{t}}\right),
\end{multline*}
or
\begin{equation}\label{sfdgdnhgfXX}
\P\big\{\DPtimeR{u,c}\leqslant t\big\}=\begin{cases}
F\big(t;\muIG,\lambdaIG,-\tfrac{1}{2}\big)\;
\big|_{\muIG=\frac{u}{\difP-c},\lambdaIG=\frac{u^2}{\sigma^2}},&0<c\leqslant\cS=\difP,
\\[8pt]
\exp\left\{-\dfrac{2\lambdaIG}{\HmuIG}\right\}\,F\big(t;\HmuIG,\lambdaIG,\tfrac{1}{2}\big)\;
\big|_{\HmuIG=-\frac{u}{\difP-c},\lambdaIG=\frac{u^2}{\sigma^2}},&c>\cS=\difP,
\end{cases}
\end{equation}
where
\begin{equation}\label{21w4e5hy}
F\left(x;\muIG,\lambdaIG,-\tfrac{1}{2}\right)
=\UGauss{0}{1}\left(\sqrt{\frac{\lambdaIG}{x}}\left(\frac{x}{\muIG}-1\right)\right)
+\exp\bigg\{\frac{2\lambdaIG}{\muIG}\bigg\}\,\UGauss{0}{1}
\left(-\sqrt{\frac{\lambdaIG}{x}}\left(\frac{x}{\muIG}+1\right)\right)
\end{equation}
denotes c.d.f. of (proper) inverse Gaussian distribution with parameters
$\muIG>0$ and $\lambdaIG>0$. It is ``inverse'' in that sense that while the
Gaussian distribution refers to a Brownian motion's position at a fixed time,
the inverse Gaussian distribution refers to the time a diffusion process takes
to reach a fixed level.

%%%%%%%%%%%%%%%%%%%%%%%%%%%%%%%%%%%%%%%%%%%%%%%%%%%%%%%%%%%%%%%%%%%%%%%%%%%%%%%%%%%%%%%%
\subsection{Level crossing by a compound renewal process}\label{sdrthy5trj}
%%%%%%%%%%%%%%%%%%%%%%%%%%%%%%%%%%%%%%%%%%%%%%%%%%%%%%%%%%%%%%%%%%%%%%%%%%%%%%%%%%%%%%%%

Let us denote by $f_{\T{1}}(t)$ and $f_{T}(t)$ the probability density
functions (p.d.f.) of a positive random variable $\T{1}$ and of a set of
positive random variables $\T{i}\overset{d}{=}T$, $i=2,3,\dots$, all
distributed identically. Introducing compound renewal process, the random
variable $\T{1}$ is the time interval between starting time zero and time of
the first renewal, and the random variables $\T{i}$ are the inter-renewal
times. The distribution of $\T{1}$ may be different from the distribution of
$T$.

By $f_{Y}(t)$ we denote p.d.f. of positive random variables
$\Y{i}\overset{d}{=}Y$, $i=1,2,\dots$, all distributed identically. The random
variables $\Y{i}$ are sizes of jumps which occur only in the moments of
renewals. Throughout the entire presentation, p.d.f. $f_{T}(y)$ and $f_{Y}(y)$
are assumed bounded from above by a finite constant.

Having assumed that $\T{1}$, i.i.d. $\T{i}\overset{d}{=}T$, $i=2,3,\dots$,
i.i.d. $\Y{i}\overset{d}{=}Y$, $i=1,2,\dots$, are all mutually independent, we
are within the renewal model, where compound renewal process with time
$s\geqslant 0$ is $\homV{s}=\sum_{i=1}^{\homN{s}}\Y{i}$, or $0$, if
$\homN{s}=0$ (or $\T{1}>s$), where
$\homN{s}=\max\left\{n>0:\sum_{i=1}^{n}\T{i}\leqslant s\right\}$, or $0$, if
$\T{1}>s$. The random variable
\begin{equation*}
\RtimeR{u,c}=\inf\left\{s>0:\homV{s}-cs>u\right\},
\end{equation*}
or $+\infty$, as $\homV{s}-cs\leqslant u$ for all $s>0$, is the first passage
time to level $u>0$ of the process $\homV{s}-cs$, $s>0$.

Let us denote by $\P\{\RtimeR{u,c}\leqslant t\mid\T{1}=v\}$ the distribution of
$\RtimeR{u,c}$ conditioned by $\T{1}=v$. It is easily seen that for $0<v<t$
\begin{multline}\label{ewrktulertye}
\P\{\RtimeR{u,c}\leqslant t\}=\int_{0}^{t}\P\{u+cv-Y<0\}\,f_{\T{1}}(v)dv
\\[-2pt]
+\int_{0}^{t}\P\{v<\RtimeR{u,c}\leqslant t\mid\T{1}=v\}\,f_{\T{1}}(v)dv.
\end{multline}
Furthermore (see Theorem 2.1 in \cite{[M1]} and references therein),
\begin{equation}\label{dfgbehredf}
\P\big\{v<\RtimeR{u,c}\leqslant
t\mid\T{1}=v\big\}=\int_{v}^{t}\dfrac{u+cv}{u+cz}\sum_{n=1}^{\infty}
\P\big\{M(u+cz)=n\big\}\,f_{T}^{*n}(z-v)dz,
\end{equation}
where $M(x)=\inf\big\{k\geqslant 1:\sum_{i=1}^{k}Y_{i}>x\big\}-1$.

%%%%%%%%%%%%%%%%%%%%%%%%%%%%%%%%%%%%%%%%%%%%%%%%%%%%%%%%%%%%%%%%%%%%%%%%%%%%%%%%%%%%%%%%
\section{Explicit solutions in level crossing}\label{gfyjrtjkrt}
%%%%%%%%%%%%%%%%%%%%%%%%%%%%%%%%%%%%%%%%%%%%%%%%%%%%%%%%%%%%%%%%%%%%%%%%%%%%%%%%%%%%%%%%

Similarly to diffusion set-up, in the compound renewal framework with $T$, $Y$
exponential with parameters $\paramT$, $\paramY$, the distribution of
$\RtimeR{u,c}$ may be written explicitly. Assuming that $\T{1}\overset{d}{=}T$
and writing
$\BesselI{n}(x)=\sum_{k=0}^{\infty}\frac{1}{k!\,(n+k)!}\,(\frac{x}{2})^{n+2k}$
for the modified Bessel function of the first kind of order $n$ (see, e.g.,
\cite{[Abramowitz Stegun 1972]}), we have the following equivalent formulas
\eqref{fghnfgmngm}--\eqref{sadfgbndfn} which all may be derived\footnote{See,
in particular, Theorem 2.2 in \cite{[M2]}. The proof of equivalence of Type~I
and Type~II formulas applies Lommel's formula.} from \eqref{ewrktulertye} and
\eqref{dfgbehredf}: Type~I formula
\begin{multline}\label{fghnfgmngm}
\P\big\{\RPtimeR{u,c}\leqslant
t\big\}=e^{-u\paramY}\paramT\int_{0}^{t}e^{-(\paramY
c+\paramT)x}\Big(\BesselI{0}\big(2\sqrt{\paramT\paramY x(cx+u)}\big)
\\[-2pt]
-\frac{cx}{cx+u}\BesselI{2}\big(2\sqrt{\paramT\paramY x(cx+u)}\big)\Big)\,dx,
\end{multline}
Type~II formula
\begin{multline}\label{sarfgbery}
\P\big\{\RPtimeR{u,c}\leqslant
t\big\}=e^{-u\paramY}\frac{\sqrt{\paramT}}{\sqrt{c\paramY}}\int_{0}^{c\paramY
t}e^{-(1+\paramT/(c\paramY))x}
\\[-2pt]
\times\sum_{n=0}^{\infty}\frac{u^n}{n!}\left(\frac{\paramT\paramY}{c}\right)^{n/2}
\frac{n+1}{x}\BesselI{n+1}\big(2x\sqrt{\paramT/(c\paramY)}\,\big)\,dx,
\end{multline}
and Type~III formula
\begin{equation}\label{sadfgbndfn}
\P\big\{\RPtimeR{u,c}\leqslant
t\big\}=\P\big\{\RPtimeR{u,c}<\infty\big\}-\frac{1}{\pi}\int_0^{\pi}f(x,u,t)\,dx,
\end{equation}
where
\begin{equation*}
\P\big\{\RPtimeR{u,c}<\infty\big\}=\begin{cases}
\dfrac{\paramT}{c\paramY}\exp\{-u(c\paramY-\paramT)/c\},&\paramT/(c\paramY)<1,
\\[6pt]
1,&\paramT/(c\paramY)\geqslant 1,
\end{cases}
\end{equation*}
and
\begin{equation*}
\begin{aligned}
f(x,u,t)&=(\paramT/(c\paramY))\Big(1+\paramT/(c\paramY)-2\sqrt{\paramT/(c\paramY)}\cos
x\Big)^{-1}
\\
&\times\exp\Big\{u\paramY\big(\sqrt{\paramT/(c\paramY)}\cos x-1\big)-t\paramT
(c\paramY/\paramT)\big(1+\paramT/(c\paramY)-2\sqrt{\paramT/(c\paramY)}\cos
x\big)\Big\}
\\
&\times\Big(\cos\big(u\paramY\sqrt{\paramT/(c\paramY)}\sin x\big)
-\cos\big(u\paramY\sqrt{\paramT/(c\paramY)}\sin x+2x\big)\Big).
\end{aligned}
\end{equation*}

The equivalent Type~I--Type~III explicit formulas
\eqref{fghnfgmngm}--\eqref{sadfgbndfn}, as well as equation \eqref{sfdgdnhgfXX}
in diffusion set-up, may seem very cumbersome and non-informative for intuitive
understanding of the level crossing phenomenon. However, the formulas
\eqref{fghnfgmngm}--\eqref{sadfgbndfn} do comply with the following observation
done in \cite{[Feller 1971]}, Ch.~II, \S~7:
\begin{quotation}
{\small surprisingly many explicit solutions in diffusion theory, queuing
theory, and other applications involve Bessel functions. It is usually far from
obvious that the solutions represent probability distributions, and the
analytic theory required to derive their Laplace transforms and other relations
is rather complex. Fortunately, the distributions in question (and many more)
may be obtained by simple randomization procedures. In this way many relations
lose their accidental character, and much hard analysis can be avoided.}
\end{quotation}
To clarify this idea, let us show (see as well \cite{[Malinovskii Kosova
2014]}) that Type~II formula connects the problem of level crossing with random
walk with random displacements.

Recall that the random walk with random displacements (see, e.g., \cite{[Takacs
1967]}, Chapter 4, \S~22) is defined as follows. Suppose that a particle
performs a random walk\index{Process!random walk!with random displacements} on
the $X$-axis. Starting at the origin, in each step the particle moves either a
unit distance to the right with probability $p$ or a unit distance to the left
with probability $q$ ($p+q=1$, $0<p<1$). Suppose that the displacements of the
particle occur at random times in the time interval $(0,\infty)$. Denote by
$\nu(z)$ the number of steps taken in the interval $(0,z]$. We suppose that
$\{\nu(z),0\leqslant z<\infty\}$ is a Poisson process\index{Process!Poisson} of
density $1/p$ and that the successive displacements are independent of each
other and independent of the process $\{\nu(z),0\leqslant z<\infty\}$. Denote
by $\xi_{p}(z)$ the position of the particle at time $z$. In this case
$\{\xi_{p}(z),0\leqslant z<\infty\}$ is a stochastic process having stationary
independent increments, $\P\{\xi_{p}(0)=0\}=1$ and almost all sample functions
of $\{\xi_{p}(z),0\leqslant z<\infty\}$ are step functions having jumps of
magnitude $1$ and $-1$.

In this random walk model, for $y>0$ and $k=0,\pm 1,\pm 2,\dots$, we have (see
equality (3) in \cite{[Takacs 1967]})
\begin{equation}\label{srdfghtjnfr}
\begin{aligned}
\P\big\{\xi_{p}(y)=k\big\}&=e^{-y/p}\big(p/q\big)^{k/2}\BesselI{k}\big(2y\sqrt{q/p}\,\big)
\\
&=\frac{y}{k}\,\hv{k}(y\mid p),
\end{aligned}
\end{equation}
where $\hv{k}(y\mid
p)=\left(\dfrac{p}{q}\right)^{k/2}\dfrac{k}{y}\,e^{-y/p}\BesselI{k}\big(2y\sqrt{q/p}\,\big)$.
Direct calculation yields
\begin{equation*}
\E\xi_{p}(y)=\big(1-(q/p)\big)y,\quad\D\xi_{p}(y)=y/p.
\end{equation*}
For $y>0$ and integer $k>0$, we have (see equalities (8) and (9) in
\cite{[Takacs 1967]})
\begin{equation}\label{dwrgsdgbsd}
\begin{aligned}
\P\bigg\{\sup_{0\leqslant z\leqslant
y}\xi_{p}(z)<k\bigg\}&=\P\big\{\xi_{p}(y)<k\big\}-\big(p/q\big)^k\,\P\big\{\xi_{p}(y)<-k\big\}
\\
&=1-k(p/q)^{k/2}\int_{0}^{y}e^{-z/p}
\BesselI{k}\big(2z\sqrt{q/p}\,\big)\,\frac{dz}{z}
\\[4pt]
&=1-k\int_{0}^{y}\P\big\{\xi_{p}(z)=k\big\}\,\frac{dz}{z}=1-\int_{0}^{y}\hv{k}(z\mid
p)\,dz.
\end{aligned}
\end{equation}
Manipulating with \eqref{srdfghtjnfr} and \eqref{dwrgsdgbsd}, we have the
equality
\begin{multline*}
k(p/q)^{k/2}\int_{0}^{y}e^{-z/p}\BesselI{k}\big(2z\sqrt{q/p}\,\big)\,\frac{dz}{z}
=e^{-y/p}\sum_{i=k}^{\infty}(p/q)^{i/2}\BesselI{i}\big(2y\sqrt{q/p}\,\big)
\\[0pt]
+e^{-y/p}\sum_{i=k+1}^{\infty}(q/p)^{i/2-k}\BesselI{i}\big(2y\sqrt{q/p}\,\big),
\end{multline*}
which proof by the methods of Bessel functions is not at all simple.

Denoting by $\varsigma_{k}(p)$ the first hitting time of the point $k$ by the
random walk with $p\in(0,1)$ and bearing in mind \eqref{dwrgsdgbsd}, we have
\begin{equation}\label{rthtyfrt}
\P\big\{\varsigma_{k}(p)\leqslant y\big\}=\P\bigg\{\sup_{0\leqslant z\leqslant
y}\xi_{p}(z)\geqslant k\bigg\}=\int_0^{y}\hv{k}(z\mid p)\,dz.
\end{equation}
Thus, the function $\hv{k}(y\mid p)$ introduced in \eqref{srdfghtjnfr} has a
clear probabilistic meaning. It is a probability density function of
$\varsigma_{k}(p)$, i.e., of the first hitting time of the point $k$ in the
model of random walk with random displacements. It is not a surprise that the
density $\hv{k}(y\mid p)$ is defective\footnote{See Chapter II, Section 7 and
Chapter XIV, Section 6 in \cite{[Feller 1971]}.} for $p<1/2$, i.e., when the
random walk drifts to the left, and proper for $p\geqslant 1/2$, i.e., when the
drift is absent or to the right.

Making the change of variables $y=x\paramT/(c\paramY)$ in Type~II formula
\eqref{sarfgbery}, we have
\begin{equation*}
\begin{aligned}
\P\big\{\RPtimeR{u,c}\leqslant
t\big\}&=e^{-u\paramY}\sum_{n=0}^{\infty}\frac{(u\paramY)^n}{n\,!}
\int_{0}^{t\paramT}\left(\frac{\paramT}{c\paramY}\right)^{(n+1)/2}\frac{n+1}{y}
e^{-y(c\paramY+\paramT)/\paramT}
\BesselI{n+1}\big(2y\sqrt{(c\paramY)/\paramT}\big)\,dy
\\
&=e^{-u\paramY}\sum_{n=0}^{\infty}\frac{(u\paramY)^n}{n\,!}
\int_{0}^{t\paramT}\hv{n+1}(y\mid p)\,\big|_{p =\paramT/(c\paramY+\paramT)}\,dy
\\
&=e^{-u\paramY}\sum_{n=0}^{\infty}\frac{(u\paramY)^n}{n\,!}\P\big\{\varsigma_{n+1}(p)\leqslant
t\paramT\big\}\;\big|_{p =\paramT/(c\paramY+\paramT)},
\end{aligned}
\end{equation*}
which makes a link between the random walk with random displacements and the
problem of level crossing evident.

It is noteworthy that under standard assumptions in the renewal model with $T$
and $Y$ exponential we have
\begin{equation*}
\P\big\{\RPtimeR{u,c}<\infty\big\}=\begin{cases}
1,&0<c\leqslant\paramT/\paramY,
\\[4pt]
(\paramT/(c\paramY))\exp\big\{-u\paramY\,\big(1-\paramT/(c\paramY)\big)\big\},&c>\paramT/\paramY,
\end{cases}
\end{equation*}
which follows from
\begin{equation*}
\int_{0}^{\infty}\hv{n+1}(y\mid p)\,dy=\begin{cases}
1,&p\geqslant 1/2,\\
\big(p/q\big)^{n+1},&p<1/2
\end{cases}
\end{equation*}
and
$e^{-x}\sum_{n=0}^{\infty}\frac{x^n}{n\,!}(p/q)^{n+1}=(p/q)\exp\big\{-x\,\big(1-(p/q)\big)\big\}$.

%%%%%%%%%%%%%%%%%%%%%%%%%%%%%%%%%%%%%%%%%%%%%%%%%%%%%%%%%%%%%%%%%%%%%%%%%%%%%%%%%%%%%%%%
\section{Approximations in level crossing}\label{drtfyjukgul}
%%%%%%%%%%%%%%%%%%%%%%%%%%%%%%%%%%%%%%%%%%%%%%%%%%%%%%%%%%%%%%%%%%%%%%%%%%%%%%%%%%%%%%%%

Among many types of approximations for $\P\{\RtimeR{u,c}\leqslant t\}$, the
most significant\footnote{Many other approximations, including Arfwedson's
saddlepoint approximation (see \cite{[Asmussen Albrecher 2010]}, p.~133) are
set aside, but may be also discussed in the similar way.} are ``normal'' and
``diffusion'' which refinement is ``corrected diffusion'' (see
\cite{[Asmussen]}, pp.~37 and~42, and \cite{[Asmussen Albrecher 2010]}). The
former name emphasizes the structure of the approximation, while the latter
name emphasizes the method of its construction.

%%%%%%%%%%%%%%%%%%%%%%%%%%%%%%%%%%%%%%%%%%%%%%%%%%%%%%%%%%%%%%%%%%%%%%%%%%%%%%%%%%%%%%%%
\subsection{``Normal''}\label{sdrtfjdf}
%%%%%%%%%%%%%%%%%%%%%%%%%%%%%%%%%%%%%%%%%%%%%%%%%%%%%%%%%%%%%%%%%%%%%%%%%%%%%%%%%%%%%%%%

Assuming that $\T{1}\overset{d}{=}T$ and writing $X=Y-cT$, let us introduce
\begin{equation}\label{dtrhtukXX}
\mominus=\E{T}/\E{X},\quad\Dominus^2=\E(X\,\E{T}-T\,\E{X})^2/(\E{X})^3.
\end{equation}
For $0<c<\cS=\E{Y}/\E{T}$, or for\; $\E{X}=\E{Y}-c\,\E{T}>0$, we have
$\mominus>0$.

\begin{theorem}\label{sdtghjrtyjkt}
In the renewal model with\, $0<c<\cS$, we assume that $0<\Dominus^2<\infty$.
Then
\begin{equation}\label{srdutuktX}
\sup_{t>0}\,\big|\;\P\{\RtimeR{u,c}\leqslant t\}-\UGauss{\mominus u}{\Dominus^2
u}(t)\;\big|\;=\overline{o}(1),\quad u\to\infty.
\end{equation}
If, in addition, $\E(Y^3)<\infty$ and $\E(T^3)<\infty$, then this supremum is
$\underline{O}(u^{-1/2})$, as $u\to\infty$.
\end{theorem}

Let us assume that a positive solution $\adjustL$ of the equation
$\E{e^{rX}}=1$ exists. This assumption entails that the tail of the random
variable $Y$ must decrease exponentially fast\footnote{This clue assumption
does not allow extension of this approach to $Y$ with power tail.}. Let us
introduce the associated random variables $\Ass{X}$, $\Ass{T}$ which joint
distribution is $F_{\Ass{X}\Ass{T}}(dz,dw)=e^{\adjustL z}F_{XT}(dz,dw)$, and
\begin{equation}\label{drfgdthdh}
\begin{gathered}
\moplus=\E\Ass{T}/\E\Ass{X},
\quad\Doplus^2=\E(\Ass{X}\E\Ass{T}-\Ass{T}\E\Ass{X}\big)^2/(\E\Ass{X})^3,
\\[4pt]
C=\frac{1}{\adjustL\,\E\Ass{X}}\exp\left\{-\sum^{\infty}_{n=1}\frac{1}{n}\,
\P\{\Sum{n}>0\}-\sum^{\infty}_{n=1}\frac{1}{n}\,\P\{\Ass{S}_{n}\leqslant
0\}\right\},
\end{gathered}
\end{equation}
where $\Ass{S}_{n}=\sum_{i=1}^n\Ass{X}_{i}$, $\Sum{n}=\sum_{i=1}^{n}\X{i}$. For
$c>\cS$, we have $\E\Ass{X}>0$ and $\moplus>0$.

\begin{theorem}\label{sdgehtdrth}
In the renewal model with\, $c>\cS$, we assume that $0<\Doplus^2<\infty$. Then
for $C>0$ defined in \eqref{drfgdthdh}
\begin{equation}\label{wery5u65X}
\sup_{t>0}\,\big|\;e^{\adjustL u}\,\P\{\RtimeR{u,c}\leqslant
t\}-C\,\UGauss{\moplus u}{\Doplus^2 u} (t)\;\big|\;=\overline{o}(1),\quad
u\to\infty.
\end{equation}
If, in addition, $\E(Y^3)<\infty$ and $\E(T^3)<\infty$, then this supremum is
$\underline{O}(u^{-1/2})$, as $u\to\infty$.
\end{theorem}

Developing the classical results by Lundberg and Cram{\'e}r, the approximation
\eqref{wery5u65X} was first obtained by Segerdahl \cite{[Segerdahl]}. There
exist different approaches to the proof (see, e.g., Siegmund \cite{[Siegmund
one]}, von Bahr \cite{[Bahr]}). Refinements in terms of Edgeworth expansions
were first discussed\footnote{The formula (4.7) in \cite{[Asmussen]} which
``invokes the higher cumulants of $\tau$ [i.e., time to ruin --
\emph{V.M.}]$\dots$ and uses an Edgeworth expansion $\dots$ to produce
correction terms'', and which was put forth as ``first-order correction, which
suggests a formula of the type (4.7)'' differs from rigorously constructed
Edgeworth expansion (see \cite{[Malinovskii 1994]}); the method of proof in
\cite{[Malinovskii 1994]} is essentially similar to the method in \cite{[M1]},
\cite{[M2]}.} in Asmussen \cite{[Asmussen]} and then proved in Malinovskii
\cite{[Malinovskii 1994]}.

For $Y$ and $T$ exponential with parameters $\paramY>0$ and $\paramT>0$
respectively, easy calculation yields $\cS=\paramT/\paramY$,
\begin{equation}\label{srdthgtjrf}
\begin{aligned}
C&=\paramT/(c\paramY),&\adjustL&=\paramY(1-\paramT/(c\paramY)),
\\
\mominus&=-\frac{1}{c(1-\paramT/(c\paramY))},&
\Dominus^2&=-\frac{2(\paramT/(c\paramY))^2+(1-\paramT/(c\paramY))^2}{\paramT
c(1-\paramT/(c\paramY))^3},
\\
\moplus&=\frac{(\paramT/(c\paramY)}{c(1-\paramT/(c\paramY))},&
\Doplus^2&=\frac{2(\paramT/(c\paramY))}{c^2\paramY(1-\paramT/(c\paramY))^3}.
\end{aligned}
\end{equation}
The following result follows straightforwardly from Theorems~\ref{sdtghjrtyjkt}
and \ref{sdgehtdrth}.

\begin{figure}[t]
\includegraphics[scale=.75]{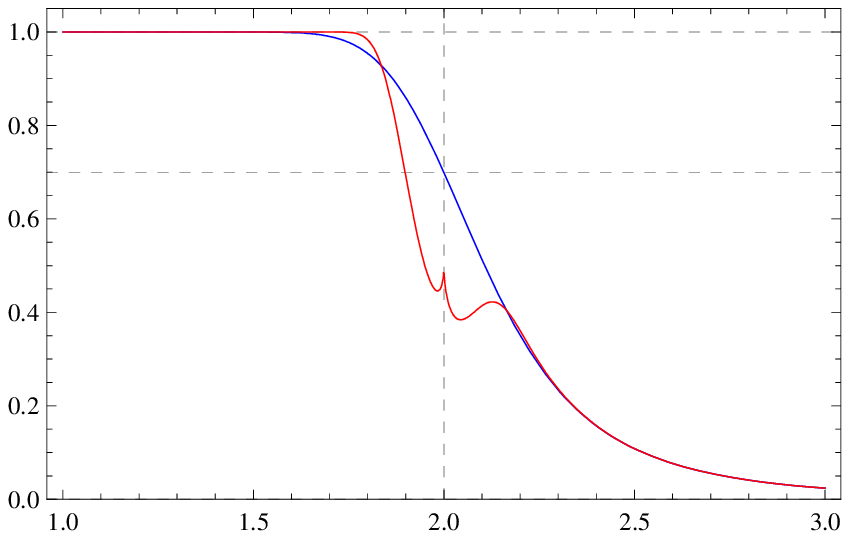}
\caption{\small Graphs ($X$-axis is $c$) of $\P\big\{\RPtimeR{u,c}\leqslant
t\big\}$ (blue line) calculated numerically using Type~I exact formula
\eqref{fghnfgmngm} and of the approximations of Theorems~\ref{etrjykih} and
\ref{sdgehtdrth} (red line) when $Y$ and $T$ are exponential with parameters
$\paramY=1$ and $\paramT=2$, $t=200$, $u=10$. Horizontal line:
$\P\big\{\RPtimeR{u,\cS}\leqslant t\big\}=0.699$.}\label{dthtjtjkt}
\end{figure}

\begin{figure}[t]
\includegraphics[scale=.75]{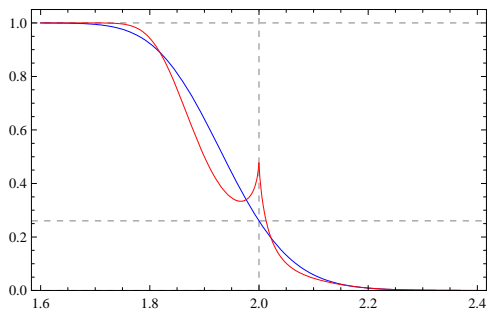}
\caption{\small Graphs as in Fig.~\ref{dthtjtjkt} for $t=500$, $u=50$.
Horizontal line is $0.26$.}\label{erthrrjnr}
\end{figure}

\begin{corollary}\label{etrjykih}
In the renewal model with $Y$ and $T$ exponential with parameters $\paramY>0$
and $\paramT>0$, for $0<c<\cS=\paramT/\paramY$ we have
\begin{equation*}
\sup_{t>0}\,\big|\;\P\{\RPtimeR{u,c}\leqslant t\}-\UGauss{\mominus
u}{\Dominus^2 u}(t)\;\big|\;=\overline{o}(1),\quad u\to\infty,
\end{equation*}
and for $c>\cS$ we have
\begin{equation*}
\sup_{t>0}\,\big|\;e^{\adjustL u}\,\P\{\RPtimeR{u,c}\leqslant
t\}-C\,\UGauss{\moplus u}{\Doplus^2 u} (t)\;\big|\;=\overline{o}(1),\quad
u\to\infty,
\end{equation*}
with $\mominus>0$, $\Dominus^2>0$, $0<C<1$, $\adjustL>0$, $\moplus>0$,
$\Doplus^2>0$ defined in \eqref{srdthgtjrf}.
\end{corollary}

The case $c=\cS$ is excluded from consideration in both
Theorems~\ref{sdtghjrtyjkt} and~\ref{sdgehtdrth}. Poor performance of the
``normal'' approximation around $\cS$ is illustrated in Figs.~\ref{dthtjtjkt}
and~\ref{erthrrjnr}.

In applications, the case $c=\cS$ is known under different names. Related is
the term ``heavy traffic''. It comes from queueing theory, but has an obvious
interpretation also in risk theory: on the average, the premiums $cs$ within
time $s>0$ exceed only slightly the expected claims $(\E{Y}/\E{T})s$ within the
same time $s>0$. That is, heavy traffic conditions mean that the safety loading
$\tau=(\E{T}/\E{Y})c-1$ is positive but small.

The approximation of Theorem~\ref{sdgehtdrth} when positive $\tau$ depends on
$u$ and tends to zero, as $u\to\infty$, was investigated in \cite{[Malinovskii
2000]}. In \cite{[Malinovskii Kosova 2014]}, further insight into performance
of Theorems~\ref{sdtghjrtyjkt} and~\ref{sdgehtdrth}, with rationale of its
poorness in the vicinity of $\cS$, may be found.

%%%%%%%%%%%%%%%%%%%%%%%%%%%%%%%%%%%%%%%%%%%%%%%%%%%%%%%%%%%%%%%%%%%%%%%%%%%%%%%%%%%%%%%%
\subsection{``Diffusion''}\label{sdrteergh}
%%%%%%%%%%%%%%%%%%%%%%%%%%%%%%%%%%%%%%%%%%%%%%%%%%%%%%%%%%%%%%%%%%%%%%%%%%%%%%%%%%%%%%%%

The idea behind the ``diffusion'' approximation is to first approximate the
claim surplus process by a Brownian motion with drift by matching the two first
moments, and next to note that such an approximation in particular implies that
the first passage probabilities are close\footnote{This and three following
paragraphs are quotations from \cite{[Asmussen Albrecher 2010]}.}
(\cite{[Asmussen Albrecher 2010]}, p. 136).

The idea behind the simple ``diffusion'' approximation is to replace the risk
process by a Brownian motion (by fitting the two first moments) and use the
Brownian first passage probabilities as approximation for the ruin
probabilities. Since Brownian motion is skip-free, this idea ignores (among
other things) the presence of the overshoot, which we have seen to play an
important role for example for the Cram{\'e}r--Lundberg approximation. The
objective of the corrected ``diffusion'' approximation is to take this and
other deficits into consideration (\cite{[Asmussen Albrecher 2010]}, p. 139).

Diffusion approximations of random walks via Donsker's theorem\footnote{Called
also Donsker--Prohorov's Invariance Principle, with credits to
\cite{[Prokhorov]}. --- \emph{V.M.}} is a classical topic of probability
theory. See for example Billingsley \cite{[Billingsley]}. The first application
in risk theory is Iglehart \cite{[Iglehart]} and two further standard
references in the area are Grandell \cite{[Grandell one]}, \cite{[Grandell
two]}. All material of this section can be found in these references. For
claims with infinite variance, Furrer, Michna, Weron \cite{[Furrer Michna
Weron]} suggested an approximation by a stable L{\'e}vy process rather than a
Brownian motion. Further relevant references in this direction are Furrer
\cite{[Furrer]} Boxma, Cohen \cite{[Boxma Cohen]} and Whitt \cite{[Whitt]}
(\cite{[Asmussen Albrecher 2010]}, p. 139)

A ``corrected diffusion'' approximations were introduced by Siegmund
\cite{[Siegmund]} in a discrete random walk setting, with the translation to
risk processes being carried out by Asmussen \cite{[Asmussen]}; this case is in
part simpler than the general random walk case because the ladder height
distribution can be found explicitly which avoids the numerical integration
involving characteristic functions which was used in \cite{[Siegmund]} to
determine the constants. In Siegmund's book \cite{[Siegmund book]} the approach
to the finite horizon case is in part different and uses local central limit
theorems. The adaptation to risk theory has not been carried out. The
``corrected diffusion'' approximation was extended to the renewal model in
Asmussen, H{\o}jgaard \cite{[Asmussen Hojgaard]}, and to the Markov-modulated
model of Chapter VII in Asmussen \cite{[Asmussen two]}; Fuh \cite{[Fuh]}
considers the closely related case of discrete time Markov additive processes.
Hogan \cite{[Hogan]} considered a variant of the ``corrected diffusion''
approximation which does not require exponential moments. His ideas were
adapted by Asmussen, Binswanger \cite{[Asmussen Binswanger]} to derive
approximations for the infinite horizon ruin probability when claims are
heavy-tailed; the analogous analysis of finite horizon ruin probabilities has
not been carried out and seems non-trivial. For ``corrected diffusion''
approximations with higher-order terms, see Blanchet, Glynn \cite{[Blanchet
Glynn]} their results also cover some heavy-tailed cases (\cite{[Asmussen
Albrecher 2010]}, p. 145).

Assuming that $\T{1}\overset{d}{=}T$, we focus on the case of $T$ and $Y$
exponential\footnote{It is noteworthy that ``diffusion'' and ``corrected
diffusion'' approximations are identical in this case.} with parameters
$\paramT$ and $\paramY$, where the distribution of $\RtimeR{u,c}$ may be
written explicitly. We follow the idea to use the approximation
\eqref{sfdgdnhgfXX} and to match the two first moments. We note that for
$R_{s}=(\difP-c)s+\sigma\W{s}$ we have $\E{R_{s}}=(\difP-c)s$,
$\D{R_{s}}=\sigma^2s$ and for $R_{s}=\sum_{i=1}^{\homN{s}}\Y{i}-cs$ we have
$\E{R_{s}}=(\paramT/\paramY-c)\,s$, $\D{R_{s}}=(2\paramT/\paramY^2)s$. Matching
these two first moments, means taking $\difP=\paramT/\paramY$ and
$\sigma^2=2\paramT/\paramY^2$ in \eqref{sfdgdnhgfXX}. It
yields\footnote{Notation $f(x;z)\,|_{z=a}$ means $f(x;a)$. Notation
$f(x)\,|_{x=a}^{b}$ means the difference $f(b)-f(a)$.}
\begin{equation}\label{asrdhrtj}
\P\big\{\RPtimeR{u,c}\leqslant t\big\}\approx\begin{cases}
F\big(t;\muIG,\lambdaIG,-\tfrac{1}{2}\big)\; \big|_{\muIG=\frac{\paramY
u}{\paramT-c\paramY},\lambdaIG=\frac{\paramY^2
u^2}{2\paramT}},&0<c\leqslant\cS,
\\[8pt]
\exp\bigg\{-\dfrac{2\lambdaIG}{\HmuIG}\bigg\}\,F\big(t;\HmuIG,\lambdaIG,\tfrac{1}{2}\big)\;
\big|_{\HmuIG=-\frac{\paramY
u}{\paramT-c\paramY},\lambdaIG=\frac{u^2}{\sigma^2}},&c>\cS,
\end{cases}
\end{equation}
where $\cS=\paramT/\paramY$, and
\begin{equation}\label{qsrfghbfds}
\muIG=\frac{u}{\difP-c}=\frac{u}{(\paramT/\paramY)-c},\quad
\lambdaIG=\frac{u^2}{\sigma^2}=\frac{\paramY^2 u^2}{2\paramT},\quad
\HmuIG=-\frac{u}{\difP-c}=-\frac{u}{(\paramT/\paramY)-c}.
\end{equation}

The approximation \eqref{asrdhrtj} illustrated in Fig.~\ref{edtfyjhm} below is
further discussed in Section~\ref{tryjukhg}.

%%%%%%%%%%%%%%%%%%%%%%%%%%%%%%%%%%%%%%%%%%%%%%%%%%%%%%%%%%%%%%%%%%%%%%%%%%%%%%%%%%%%%%%%
\section{New ``inverse Gaussian'' approximation in level crossing}\label{dsrtttgf}
%%%%%%%%%%%%%%%%%%%%%%%%%%%%%%%%%%%%%%%%%%%%%%%%%%%%%%%%%%%%%%%%%%%%%%%%%%%%%%%%%%%%%%%%

\begin{figure}[t]
\includegraphics[scale=0.8]{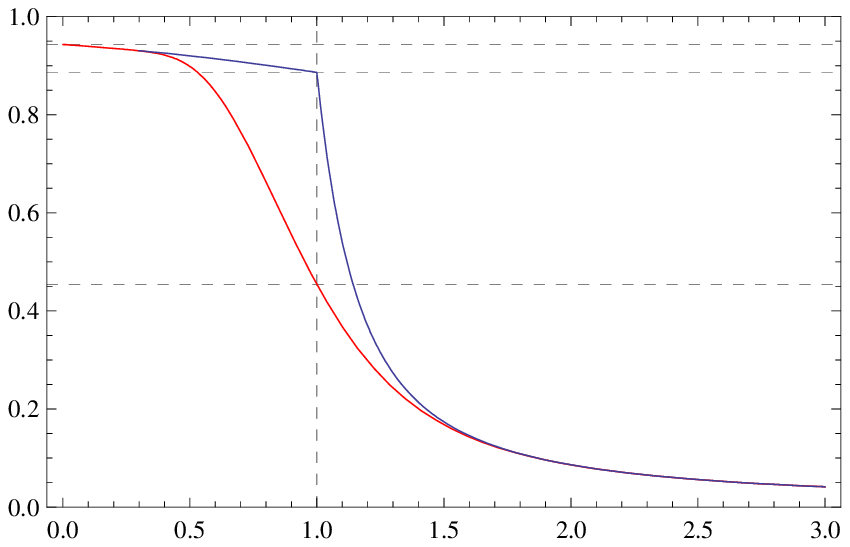}
\caption{\small Graphs ($X$-axis is $c$) of the functions
$\AInt{M}{u,c}{\infty\mid 0}$ (blue line) and $\AInt{M}{u,c}{t\mid 0}$ (red
line), as $M=1$, $D^2=6$, $t=100$, $u=15$. Horizontal lines are
$\AInt{M}{u,0}{\infty\mid 0}=0.943$, $\AInt{M}{u,\cS}{\infty\mid 0}=0.886$, and
$\AInt{M}{u,\cS}{t\mid 0}=0.454$, where $\cS=\E{Y}/\E{T}$.}\label{sdrtgerher}
\end{figure}

Being in the renewal framework of Section \ref{sdrthy5trj}, let us write
$M={\E{T}}/{\E{Y}}$, $D^2=((\E{T})^2\D{Y}+(\E{Y})^2\D{T})/(\E{Y})^3$, and
\begin{equation}\label{wqerthrrthj}
\AInt{M}{u,c}{t\mid v}=\int_{0}^{\frac{c(t-v)}{u+cv}}\frac{1}{1+y}
\,\Ugauss{cM(1+y)}{\frac{c^2D^2(1+y)}{u+cv}}(y)dy,
\end{equation}
which is expressed in terms of c.d.f. of inverse Gaussian distribution
\eqref{21w4e5hy} as
\begin{equation}\label{sdfghdtnYY}
\AInt{M}{u,c}{t\mid v}=\begin{cases}
F\big(x;\muIG,\lambdaIG,-\tfrac{1}{2}\big)\;\Big|_{x=1}^{\frac{c(t-v)}{cv+u}+1},&0<c\leqslant\cS,
\\[10pt]
\exp\Bigg\{-2\,\dfrac{\lambdaIG}{\HmuIG}\Bigg\}
F\big(x;\HmuIG,\lambdaIG,-\tfrac{1}{2}\big)\;\Big|_{x=1}^{\frac{c(t-v)}{cv+u}+1},&
c\geqslant\cS,
\end{cases}
\end{equation}
where $\cS=\frac{1}{M}$, $\lambdaIG=\frac{u}{c^2D^2}>0$,
$\muIG=\frac{1}{1-cM}$, and $\HmuIG=\frac{1}{cM-1}$. Plainly, $\muIG>0$ for
$0<c<\cS$ and $\HmuIG>0$ for $c>\cS$.

%%%%%%%%%%%%%%%%%%%%%%%%%%%%%%%%%%%%%%%%%%%%%%%%%%%%%%%%%%%%%%%%%%%%%%%%%%%%%%%%%%%%%%%%
\subsection{Approximation for conditional distribution}\label{ertuyuyryt}
%%%%%%%%%%%%%%%%%%%%%%%%%%%%%%%%%%%%%%%%%%%%%%%%%%%%%%%%%%%%%%%%%%%%%%%%%%%%%%%%%%%%%%%%

The following result is formulated and proved as Theorem 1.1 in \cite{[M1]}.

\begin{theorem}\label{asrghjtrn}
\label{srdthjrf} In the renewal model, let p.d.f. $f_{T}(y)$ and $f_{Y}(y)$ be
bounded from above by a finite constant, $D^2>0$, $\E({T}^3)<\infty$,
$\E({Y}^3)<\infty$. Then for fixed $c>0$ and $0<v<t$ we have
\begin{equation}\label{werhrjet}
\sup_{t>v}\;\big|\;\P\big\{v<\RtimeR{u,c}\leqslant
t\mid\T{1}=v\big\}-\AInt{M}{u,c}{t\mid v}\;\big|\;
=\underline{O}\,\left(\frac{\ln(u+cv)}{u+cv}\right),
\end{equation}
as $u+cv\to\infty$.
\end{theorem}

Since $\AInt{M}{u,c}{t\mid v}$ in \eqref{sdfghdtnYY} is defined by means of
c.d.f. of inverse Gaussian distribution, let us call\footnote{The ``diffusion''
approximation also involves the inverse Gaussian distribution. Someone will not
concur with the proposed name, but to us it seems sensible.} the approximation
of Theorem~\ref{asrghjtrn}, as well as its corollaries put forth below,
``inverse Gaussian''.

Bearing in mind easy equality
\begin{multline*}
\P\{M(u+cv+cy)=n\}=\P\bigg\{\sum_{i=1}^{n}\Y{i}\leqslant
u+cv+cy<\sum_{i=1}^{n+1}\Y{i}\bigg\}
\\
=\int_{0}^{u+cv+cy}f^{*n}_{Y}(u+cv+cy-z)\P\{\Y{n+1}>z\}dz,
\end{multline*}
the idea behind the proof of Theorem~\ref{asrghjtrn} consists in application of
the identity \eqref{dfgbehredf}. After the change of variables $y=z-v$, it
writes as
\begin{equation*}
\begin{aligned}
\P\big\{v<\RtimeR{u,c}\leqslant
t\mid\T{1}=v\big\}&=\int_{0}^{t-v}\dfrac{u+cv}{u+cv+cy}\,\,
\p_{\big\{\sum_{i=2}^{M(u+cv+cy)+1}\T{i}\big\}}(y)dy
\\
&=\sum_{n=1}^{\infty}\int_{0}^{t-v}\frac{u+cv}{u+cv+cy}
\int_{0}^{u+cv+cy}\P\left\{\Y{n+1}>z\right\}
\\[2pt]
&\hskip 90pt\times f_{Y}^{*n}(u+cv+cy-z)f_{T}^{*n}(y)dydz,
\end{aligned}
\end{equation*}
and non-uniform central limit theorem applied to approximate $n$-fold
convolutions $f_{Y}^{*n}$, $f_{T}^{*n}$ yields the approximating and residual
terms, both investigated in \cite{[M1]} by direct analytical method. This
method is a considerable development of the method used previously in
\cite{[Malinovskii 1994]}, \cite{[Malinovskii 1996]}, and \cite{[Malinovskii
2000]}. It applies identical transformations, approximations of integral sums
by the corresponding integrals, evaluation of elliptic integrals of the third
kind and of sums related to zeta-functions.

%%%%%%%%%%%%%%%%%%%%%%%%%%%%%%%%%%%%%%%%%%%%%%%%%%%%%%%%%%%%%%%%%%%%%%%%%%%%%%%%%%%%%%%%
\subsection{Approximation for non-conditional distribution}\label{dfyjtfg}
%%%%%%%%%%%%%%%%%%%%%%%%%%%%%%%%%%%%%%%%%%%%%%%%%%%%%%%%%%%%%%%%%%%%%%%%%%%%%%%%%%%%%%%%

Though approximation for conditional distribution is more convenient to prove
because in this case the identity \eqref{dfgbehredf} is a convenient structure,
the traditional object of interest is the unconditional distribution
$\P\big\{\RtimeR{u,c}\leqslant t\big\}$. Accordingly, let us formulate some
corollaries of Theorem~\ref{asrghjtrn}.

\begin{corollary}\label{erty5trjrk}
Under conditions of Theorem \ref{srdthjrf}, with p.d.f. $f_{\T{1}}(y)$ bounded
from above by a finite constant and with $\E\T{1}<\infty$, we have
\begin{equation*}
\sup_{t>0}\,\Big|\,\P\big\{\RtimeR{u,c}\leqslant
t\big\}-\int_{0}^{t}\P\{Y>u+cv\}f_{\T{1}}(v)dv -\int_{0}^{t}\AInt{M}{u,c}{t\mid
v}\,f_{\T{1}}(v)dv\,\Big|\,=\underline{O}\,\left(\frac{\ln u}{u}\right),
\end{equation*}
as $u\to\infty$.
\end{corollary}

Having information about ``regularity'' of $\T{1}$, the insight into the term
$\int_{0}^{t}\P\{Y>u+cv\}f_{\T{1}}(v)dv$ is easy. To be particular, application
of the Markov's inequality yields $\P\{Y>u+cv\}\leqslant\E{Y^3}/(u+cv)^3$, and
$\int_{0}^{t}(u+cv)^{-3}f_{\T{1}}(v)dv=\underline{O}\,(u^{-3})$, as
$u\to\infty$. This order of magnitude is less than $\underline{O}\,(\ln u/u)$,
and this term may be omitted.

Concerning the second term, we bear in mind that for each $v>0$
\begin{equation*}
\AInt{M}{u,c}{t\mid
v}=\AInt{M}{u,c}{t-v}+\underline{O}\,\bigg(\dfrac{v}{u}\bigg),\quad u\to\infty,
\end{equation*}
where $\AInt{M}{u,c}{t-v}=\AInt{M}{u,c}{t-v\mid 0}$. It yields the following
result.

\begin{corollary}\label{rwtherherj}
Under conditions of Corollary~\ref{erty5trjrk}, we have
\begin{multline}\label{sdrthyjufd}
\sup_{t>0}\,\Big|\,\P\big\{\RtimeR{u,c}\leqslant
t\big\}-\int_{0}^{t}\P\{Y>u+cv\}f_{\T{1}}(v)dv
\\[-2pt]
-\int_{0}^{t}\AInt{M}{u,c}{t-v}\,f_{\T{1}}(v)dv\,\Big|\,=\underline{O}\,\left(\frac{\ln
u}{u}\right),
\end{multline}
as $u\to\infty$.
\end{corollary}

The term $\int_{0}^{t}\AInt{M}{u,c}{t-v}\,f_{\T{1}}(v)dv$ which corresponds to
the case $\T{1}<\RtimeR{u,c}$, is a convolution. It agrees with the
probabilistic intuition about the r{\^o}le of the first interval $\T{1}$ in the
event of first level $u$ crossing: when $\T{1}$ is fixed and is equal to $v$,
the time is modified from the whole time $t$ to reduced time $t-v$.

Until now we did not assume that time $t$ is large, but deemed that it may be
small, moderate, and large. Since the approximation \eqref{sdrthyjufd} is
formulated uniformly with respect to $t>0$, the influence of $\T{1}$ can not be
eliminated for $t$ small and moderate. Though $\AInt{M}{u,c}{t-v}$ is rendered
in a relatively compact form, the integral
$\int_{0}^{t}\AInt{M}{u,c}{t-v}\,f_{\T{1}}(v)dv$ hardly can be evaluated as a
compact explicit expression, even for $\T{1}$ exponential. On the contrary, for
large $t$ and ``regular'' $\T{1}$, the first interval does not affect the
main-term approximation, i.e., for $t\to\infty$ we have
\begin{equation}\label{dsfgnhrgmn}
\int_{0}^{t}\AInt{M}{u,c}{t-v}\,f_{\T{1}}(v)dv=\AInt{M}{u,c}{t}\,(1+\overline{o}(1)).
\end{equation}

\begin{remark}\label{wsdrthyjf}
The point $\cS$ is excluded from consideration in Theorems~\ref{sdtghjrtyjkt}
and~\ref{sdgehtdrth}. In this respect, it is special. In this point, the
approximation \eqref{sdrthyjufd} for $\P\big\{\RtimeR{u,c}\leqslant t\big\}$,
being neither outstanding, nor especially good, but merely noteworthy, involves
the expression
\begin{equation*}
\AInt{M}{u,\cS}{t}=2\left(\UGauss{0}{1}\left(\dfrac{\sqrt{u}}{\cS
D}\;\right)-\UGauss{0}{1}\left(\dfrac{u}{\cS D\sqrt{\cS t+u}}\;\right)\right).
\end{equation*}
The approximation for $\P\big\{\RtimeR{u,c}\leqslant\infty\big\}$ involves the
expression
\begin{equation*}
\AInt{M}{u,\cS}{\infty}= 2\,\UGauss{0}{1}\left(\frac{M\sqrt{u}}{D}\,\right)-1.
\end{equation*}
\end{remark}

%%%%%%%%%%%%%%%%%%%%%%%%%%%%%%%%%%%%%%%%%%%%%%%%%%%%%%%%%%%%%%%%%%%%%%%%%%%%%%%%%%%%%%%%
\subsection{Focus of the approximation and related problems}\label{sdyjumt}
%%%%%%%%%%%%%%%%%%%%%%%%%%%%%%%%%%%%%%%%%%%%%%%%%%%%%%%%%%%%%%%%%%%%%%%%%%%%%%%%%%%%%%%%

The approximation \eqref{sdrthyjufd} is informative when
$\P\big\{\RtimeR{u,c}\leqslant t\big\}$ tends to a positive value, rather than
to zero, as $u\to\infty$. It is paramount, e.g., in the problem of
investigating a solution $u_{\alpha,t}(c)$ of the equation
\begin{equation}\label{dsfgnhrgmn}
\P\big\{\RtimeR{u,c}\leqslant t\big\}=\alpha,
\end{equation}
where $0<\alpha<1$. For diffusion process, $u_{\alpha,t}(c)$ is investigated in
\cite{[Malinovskii 2007]}, \cite{[Malinovskii 2009]}, and \cite{[Malinovskii
2014a]}. For compound renewal process with exponential $T$ and $Y$, it was done
in \cite{[Malinovskii 2012]}, \cite{[Malinovskii 2014b]}. Analysis of
$u_{\alpha,t}(c)$ for general compound renewal process was done in
\cite{[Malinovskii Kosova 2014]} and \cite{[M3]}. It was found that the main
term in approximation of $u_{\alpha,t}(c)$, as $t\to\infty$, is proportional to
$t^{1/2}$.

The focus in \eqref{sdrthyjufd} with the remainder term $\underline{O}\,(\ln
u/u)$ differs from the focus in \eqref{wery5u65X} with the remainder term
$\overline{o}\,(e^{-\adjustL u})$, as $u\to\infty$; the latter is highly
informative when $\P\big\{\RtimeR{u,c}\leqslant t\big\}$ tends to zero, as
$u\to\infty$, while the former may be not.

\begin{remark}[Similarity with CLT]\label{sdrtgbhtr}
For approximations \eqref{sdrthyjufd} and \eqref{wery5u65X}, the difference in
focuses is similar to situation known for the normal approximation and large
deviations in the common central limit theorem for sums of i.i.d. summands. The
former concerns itself with the asymptotic behavior around the mean value of
the sum considered, while the latter deals with the exponential decline of
remote tails of the distribution of the sum.
\end{remark}

%%%%%%%%%%%%%%%%%%%%%%%%%%%%%%%%%%%%%%%%%%%%%%%%%%%%%%%%%%%%%%%%%%%%%%%%%%%%%%%%%%%%%%%%
\subsection{Subexponential distributions}\label{drtsgjhfgk}
%%%%%%%%%%%%%%%%%%%%%%%%%%%%%%%%%%%%%%%%%%%%%%%%%%%%%%%%%%%%%%%%%%%%%%%%%%%%%%%%%%%%%%%%

Subexponential distributions are a special class of heavy-tailed distributions
prominent in applied probability. First order approximations to ruin
probabilities and waiting time distributions are by now called ``folklore found
in the relevant textbooks'', e.g., Embrechts, Kl{\"u}ppelberg and Mikosch
\cite{[Embrechts Kluppelberg Mikosch]}.

For some of such distributions (e.g., Pareto inter-renewal times and Pareto
jump sizes), performance of ``inverse Gaussian'' approximation obtained in
\cite{[M1]}, \cite{[M2]} was tested against simulation in \cite{[M3]}. As we
claimed in Section~\ref{sdyjumt} above, the interest in \cite{[M3]} was the
case when $\P\{\RtimeR{u,c}\leqslant t\}$ tends to a non-zero value, as
$u\to\infty$, rather than ``large deviations'' case, i.e., when
$\P\{\RtimeR{u,c}\leqslant t\}$ tends to zero, as $u\to\infty$.

It is noteworthy that in the paper \cite{[Asmussen Kluppelberg]} called ``Large
deviations results for subexponential tails, with applications to insurance
risk'', Asmussen and Kl{\"u}ppelberg gave asymptotic expressions for
$\P\{\RtimeR{u,c}\leqslant t\}$ when $t=t(u)$, such that $t(u)/u\to
k\in(0,\infty)$, which implies that $\P\{\RtimeR{u,c}\leqslant t\}$ is tending
to zero, as $u\to\infty$. But this case is just ``large deviations'', which
lies outside the focus of \cite{[M1]}--\cite{[M3]}.

%%%%%%%%%%%%%%%%%%%%%%%%%%%%%%%%%%%%%%%%%%%%%%%%%%%%%%%%%%%%%%%%%%%%%%%%%%%%%%%%%%%%%%%%
\subsection{Extensions of \eqref{sdrthyjufd}}\label{wsdyjyt}
%%%%%%%%%%%%%%%%%%%%%%%%%%%%%%%%%%%%%%%%%%%%%%%%%%%%%%%%%%%%%%%%%%%%%%%%%%%%%%%%%%%%%%%%

In \cite{[M2]}, the asymptotic expansions are constructed with the first
correction term given explicitly. It is done by extending the same technique as
was used in \cite{[M2]}. A concurrent objective of \cite{[M2]} was to
demonstrate the amplitude of the method. The terms in which the correction is
found, are the generalized inverse Gaussian distribution. The results of
\cite{[M1]}, \cite{[M2]} are illustrated numerically in \cite{[M3]}.

\begin{remark}\label{sdfnjfmg}
There is a difference between asymptotic expansions as a mathematical result,
and as a tool to improve numerical performance. Many examples (in a simpler
problems, e.g., in the central limit theorem for sums) show that involving a
correction term, one not always and not necessary acquire a visible improvement
in terms of proximity to the approximated function.
\end{remark}

%%%%%%%%%%%%%%%%%%%%%%%%%%%%%%%%%%%%%%%%%%%%%%%%%%%%%%%%%%%%%%%%%%%%%%%%%%%%%%%%%%%%%%%%
\subsection{The asymptotic behavior, as $t\to\infty$}\label{sfdghtrjn}
%%%%%%%%%%%%%%%%%%%%%%%%%%%%%%%%%%%%%%%%%%%%%%%%%%%%%%%%%%%%%%%%%%%%%%%%%%%%%%%%%%%%%%%%

\begin{figure}[t]
\includegraphics[scale=.75]{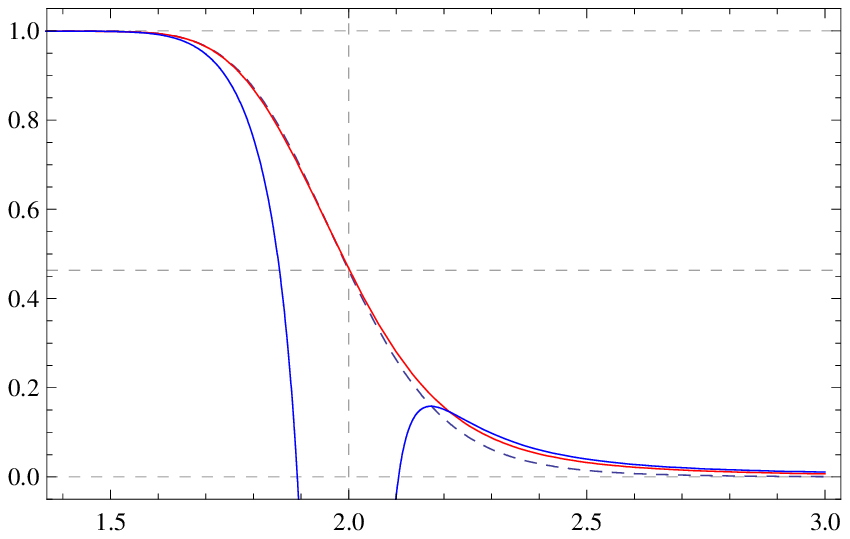}
\caption{\small Graphs ($X$-axis is $c$) of $\P\big\{\RPtimeR{u,c}\leqslant
t\big\}$ calculated numerically using Type~I exact formula \eqref{fghnfgmngm}
(dashed line), of the approximation $\AInt{M}{u,c}{t}=\AInt{M}{u,c}{t\mid 0}$
given by \eqref{sdfghdtnYY} (red line), and of approximation \eqref{asdfghjnf}
(blue line). Here $t=200$, $u=20$, $\paramY=1$, $\paramT=2$, Horizontal line:
$\P\big\{\RPtimeR{u,\cS}\leqslant t\big\}=0.463$.}\label{dfgfdgfsd}
\end{figure}

Assuming that $\homN{s}$, $s>0$, is a Poisson process and that
$c>\cS=\E{Y}/\E{T}$, Teugels \cite{[Teugels 1982]} obtained\footnote{Quoting
\cite{[Asmussen]}, ``a relation of this type is highly expected$\dots$ A
rigorous proof and an explicit evaluation of $\Do$ was, however, first provided
very recently by Teugels \cite{[Teugels 1982]}'' (\cite{[Asmussen]}, p. 51).}
an approximation of the form
\begin{equation}\label{wetru6u}
\P\big\{\RtimeR{u,c}\leqslant t\big\}\approx\P\big\{\RtimeR{u,c}<\infty\}-\Do\,
u\,t^{-3/2}e^{-\gamma_1 u-\gamma_2 t}, \quad\text{as}\ t\to\infty,
\end{equation}
where $\Do$, $\gamma_1$ and $\gamma_2$ are certain positive constants. In
\cite{[Teugels 1982]}, the renewal model was considered under the following
conditions. The i.i.d. r.v. $\T{i}\overset{d}{=}T>0$, $i=1,2,\dots$, are
exponential with parameter $\paramT$, while i.i.d. r.v.
$\Y{i}\overset{d}{=}Y>0$, $i=1,2,\dots$, are arbitrarily distributed subject to
the following condition\footnote{Called Assumption A(iv) on p. 164 of
\cite{[Teugels 1982]}.}: there exists a positive value $v$ inside the region of
convergence of $\Lambda(s)=cs-\paramT\,(1-\E e^{-sY})$, such that
$\Lambda^{\prime}(-v)=0$. This rather restrictive constrain implies that $\E
e^{-sY}$ converges as an analytical function in a neighborhood of the origin,
and the adjustment coefficient $\adjustL>0$ exists.

The proof in \cite{[Teugels 1982]} applies the Laplace transforms method. On
p.~174 of \cite{[Teugels 1982]} it is mentioned that for $Y$ exponential, one
can derive \eqref{wetru6u} straightforwardly, ``through an asymptotic expansion
of the modified Bessel function'' in Type~II formula \eqref{sarfgbery}.

Assuming that $t\to\infty$, Corollary~\ref{rwtherherj} yields ``Teugels-type''
approximation, as follows. If $u/t\to 0$ ($u\ll t$), we have
\begin{multline}\label{asdfghjnf}
\AInt{M}{u,c}{t}=\AInt{M}{u,c}{\infty}-\Do\,u\exp\left\{\gamma_1\,u\right\}
\\[-2pt]
\times\,t^{-3/2}\exp\left\{-\gamma_2\,t\right\}
\exp\left\{-\gamma_3\frac{u^2}{t}\right\}\,(1+\overline{o}(1)),\quad
t\to\infty,
\end{multline}
where $\Do=\sqrt{\dfrac{2}{\pi}}\dfrac{D}{c^{1/2}(1-cM)^2}$,
$\gamma_1=\dfrac{1-cM}{c^2D^2}$, $\gamma_2=\dfrac{(1-cM)^2}{2cD^2}$,
$\gamma_3=\dfrac{1}{2c^3D^2}$ and
\begin{equation*}
\AInt{M}{u,c}{\infty}=\begin{cases}
\UGauss{0}{1}\left(\frac{M}{D}\sqrt{u}\,\right)
-\exp\left\{\frac{2(1-cM)}{c^2D^2}\,u\right\}\,
\UGauss{0}{1}\left(-\frac{2-cM}{cD}\sqrt{u}\,\right),& 0\leqslant
c\leqslant\cS,
\\[10pt]
\exp\left\{-\frac{2(cM-1)}{c^2D^2}\,u\,\right\}\,\UGauss{0}{1}\left(-\frac{cM-2}{cD}\sqrt{u}\,\right)
-\UGauss{0}{1}\left(-\frac{M}{D}\sqrt{u}\,\right),& c>\cS,
\end{cases}
\end{equation*}

The approximation \eqref{asdfghjnf} is called ``Teugels-type'' since it differs
from the original Teugels approximation \eqref{wetru6u}. First, conditions on
the model in these two results are obviously different: ``Teugels-type'' was
obtained under much more general conditions. Second, it is assumed in
\cite{[Teugels 1982]} that $u$ is fixed and $t$ is tending to infinity, while
in Corollary~\ref{rwtherherj} it is assumed that both $u$ and $t$ are tending
to infinity. However, the structure of both approximations \eqref{wetru6u} and
\eqref{asdfghjnf} is the same, as is expected.

Indeed, quoting Teugels' Remark 8.2 on p.~174 of \cite{[Teugels 1982]}, since
``as a function of time the assumption A(iv) essentially introduces an
exponential decay'', for $W(x)=\P\{\inf_{s>0}(x+cs-\W{s})\geqslant 0\}$ and
$W(t,x)=\P\{\inf_{0<s\leqslant t}(x+cs-\W{s})\geqslant 0\}$ with $c=1$,
\begin{quotation}
{\small the extra $t^{-\gamma}$ with $\gamma=1/2$ for the compound Poisson
process and with $\gamma=3/2$ for its supremum is not surprising in view of an
approximating Wiener process. For let us consider the quantity $W(t,x)$ for a
Wiener process with drift $1$. Then (\cite{[Takacs 1967]}, p. 82) \dots
\begin{equation*}
W(t,x)-W(x)=e^{-2x}\left\{1-\UGauss{0}{1}\bigg(\frac{t-x}{\sqrt{t}}\bigg)\right\}
-\left\{1-\UGauss{0}{1}\bigg(\frac{t+x}{\sqrt{t}}\bigg)\right\}
\end{equation*}
$\dots$ which can be simplified to
\begin{equation*}
W(t,x)-W(x)\sim\sqrt{\frac{2}{\pi}}\,e^{x-t/2}x\,t^{-3/2},\quad t\to\infty.
\end{equation*}}
\end{quotation}

To have an illustration of \eqref{asdfghjnf} stated in a general renewal
framework, let us consider $T$ and $Y$ exponential with positive parameters
$\paramT$ and $\paramY$. We have
\begin{equation*}
M=\paramY\paramT^{-1},\quad D^2=2\paramY\paramT^{-2},
\end{equation*}
and
\begin{equation*}
\Do=
\frac{2}{\sqrt{\pi}}\frac{\sqrt{\paramY}}{\sqrt{c}\paramT(1-c\paramY/\paramT)^2},\quad
\gamma_1=\frac{\paramT(\paramT-c\paramY)}{2c^2\paramY},\quad\gamma_2=\frac{(\paramT-c\paramY)^2}{4
c\paramY},\quad\gamma_3=\frac{\paramT^2}{4c^3\paramY}.
\end{equation*}

Performance of Teugels-type approximation \eqref{asdfghjnf} in this case is
illustrated in Fig.~\ref{dfgfdgfsd}. It is poor in the vicinity of $\cS$; the
same is known for the original Teugels' approximation in \cite{[Teugels 1982]}.
It is no surprise since $\Do=
\frac{2}{\sqrt{\pi}}\frac{\sqrt{\paramY}}{\sqrt{c}\paramT(1-c\paramY/\paramT)^2}$
turns into infinity for $c=\cS$.

%%%%%%%%%%%%%%%%%%%%%%%%%%%%%%%%%%%%%%%%%%%%%%%%%%%%%%%%%%%%%%%%%%%%%%%%%%%%%%%%%%%%%%%%
\section{``Diffusion'' and ``inverse Gaussian'' approximations vs. exact formula when $T$, $Y$
are exponential}\label{tryjukhg}
%%%%%%%%%%%%%%%%%%%%%%%%%%%%%%%%%%%%%%%%%%%%%%%%%%%%%%%%%%%%%%%%%%%%%%%%%%%%%%%%%%%%%%%%

\begin{figure}[t]
\includegraphics[scale=0.8]{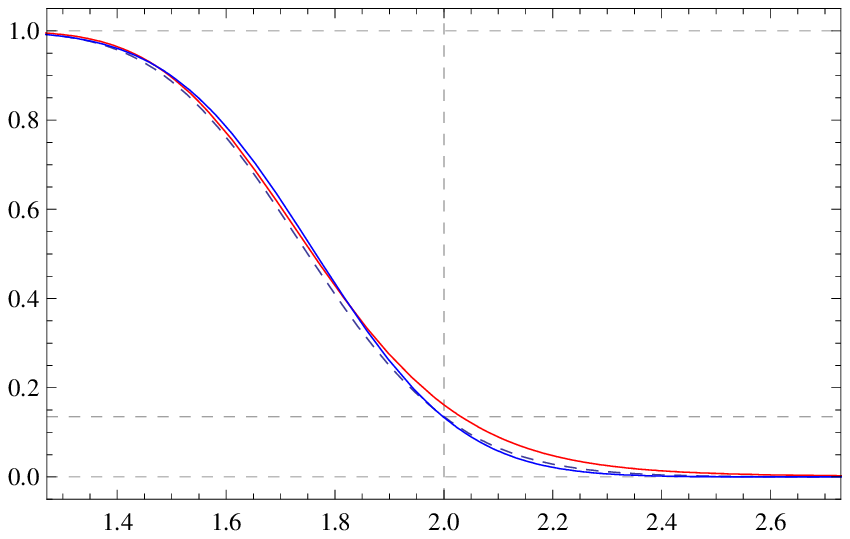}
\caption{\small Graphs ($X$-axis is $c$) of the ``diffusion'' approximation
\eqref{asrdhrtj}, \eqref{qsrfghbfds} (blue), ``inverse Gaussian'' approximation
\eqref{drtferjh}, \eqref{wsdrfgnbgf} (red) and the exact values for
$\P\big\{\RPtimeR{u,c}\leqslant t\big\}$ (dashed), as $t=100$, $u=30$,
$\paramY=1$, $\paramT=2$. Horizontal line is $\P\big\{\RPtimeR{u,\cS}\leqslant
t\big\}=0.1348$, where $\cS=\paramT/\paramY$.}\label{edtfyjhm}
\end{figure}

Bearing in mind \eqref{dsfgnhrgmn} and \eqref{sdfghdtnYY}, equation
\eqref{sdrthyjufd} for large $t>0$ yields
\begin{equation}\label{drtferjh}
\P\big\{\RPtimeR{u,c}\leqslant t\big\}\approx\begin{cases}
F\big(\frac{ct}{u}+1;\muIG,\lambdaIG,-\tfrac{1}{2}\big)-F\big(1;\muIG,\lambdaIG,-\tfrac{1}{2}\big),
&0<c\leqslant\cS,
\\[10pt]
\exp\left\{-2\,\dfrac{\lambdaIG}{\HmuIG}\right\}
\big(F\big(\frac{ct}{u}+1;\HmuIG,\lambdaIG,-\tfrac{1}{2}\big) &
\\
\hskip 120pt-F\big(1;\HmuIG,\lambdaIG,-\tfrac{1}{2}\big)\big), & c>\cS,
\end{cases}
\end{equation}
where $\cS=\frac{1}{M}$, $\muIG=\frac{1}{1-cM}$, $\lambdaIG=\frac{u}{c^2D^2}$,
$\HmuIG=\frac{1}{cM-1}$.

For $T$, $Y$ exponential with parameters $\paramT$, $\paramY$, we have
$M=\paramY\paramT^{-1}$, $D^2=2\paramY\paramT^{-2}$, and in \eqref{drtferjh} we
have $\cS=\paramT/\paramY$, and
\begin{equation}\label{wsdrfgnbgf}
\muIG=\frac{1}{1-c(\paramY/\paramT)}=\frac{(\paramT/\paramY)}{(\paramT/\paramY)-c},
\quad\lambdaIG=\frac{\paramT^2u}{2c^2\paramY}>0,\quad
\HmuIG=-\frac{1}{1-c(\paramY/\paramT)}.
\end{equation}

Having derived explicit expressions \eqref{asrdhrtj}, \eqref{qsrfghbfds} for
``diffusion'' approximation and \eqref{drtferjh}, \eqref{wsdrfgnbgf} for
``inverse Gaussian'' approximation, we are able to compare these results with
each other. Having Type~I--Type~III explicit formulas for
$\P\big\{\RPtimeR{u,c}\leqslant t\big\}$, we can further compare them with the
exact formula.

Comparing of \eqref{asrdhrtj}, \eqref{qsrfghbfds} and \eqref{drtferjh},
\eqref{wsdrfgnbgf} with each other clearly shows that the structure of
``inverse Gaussian'' approximation and the structure of ``diffusion''
approximation\footnote{Which coincides in this case with ``corrected
diffusion'' approximation.} are different, though both \eqref{drtferjh} and
\eqref{asrdhrtj} are written in terms of c.d.f. of inverse Gaussian
distribution \eqref{21w4e5hy}.

When the values of $\P\big\{\RPtimeR{u,c}\leqslant t\big\}$ are significantly
greater than zero (which happens in a region including $\cS$), performance of
both ``inverse Gaussian'' and ``diffusion'' approximations is compared to exact
values of $\P\big\{\RPtimeR{u,c}\leqslant t\big\}$ in Fig.~\ref{edtfyjhm}.

%%%%%%%%%%%%%%%%%%%%%%%%%%%%%%%%%%%%%%%%%%%%%%%%%%%%%%%%%%%%%%%%%%%%%%%%%%%%%%%%%%%%%%%%
\section{Conclusions}\label{drtsgjhfgk}
%%%%%%%%%%%%%%%%%%%%%%%%%%%%%%%%%%%%%%%%%%%%%%%%%%%%%%%%%%%%%%%%%%%%%%%%%%%%%%%%%%%%%%%%

The ``inverse Gaussian'' approximation for $\P\big\{\RtimeR{u,c}\leqslant
t\big\}$ considered in Section~\ref{dsrtttgf} is structurally
different\footnote{Though ``inverse Gaussian'', ``diffusion'', and ``corrected
diffusion'' approximations may be written in terms of c.d.f. of inverse
Gaussian distribution \eqref{21w4e5hy}.} from ``normal'', ``diffusion'', and
``corrected diffusion'' approximations (see in particular illustration in
Section~\ref{tryjukhg}).

The ``inverse Gaussian'' approximation and its refinements in terms of
Edgeworth expansions were rigorously proved by a technique based on an idea
different from the ideas lying behind the ``normal'', ``diffusion'', and
``corrected diffusion'' approximations (see Sections~\ref{sdrteergh}
and~\ref{ertuyuyryt}). The methodology is said to be creative, and cannot be
directly translated to ``standard'' techniques\footnote{Though develops the
technique in \cite{[Malinovskii 1994]} and \cite{[Malinovskii 1994]}.}; it is
not very demanding in terms of background.

The conditions which are required in ``inverse Gaussian'' approximation are
weaker, and even look minimal, in comparison with conditions in ``normal'',
``diffusion'', and ``corrected diffusion''; these substantially more
restrictive conditions seem immanent to the methods of proof for ``normal'',
``diffusion'', and ``corrected diffusion'' approximations\footnote{Quoting
again \cite{[Asmussen Albrecher 2010]} (see quotation in full in Section
\ref{sdrteergh}), note that ``the analogous analysis of finite horizon ruin
probabilities $\psi(u,T)$ has not been carried out and seems non-trivial.''}.

The ``Teugels-type'' approximation, as expected, is a corollary of the
``inverse Gaussian'' approximation; the form of the ``Teugels-type''
approximation is straightforward under very broad assumptions on the renewal
model.

The ``inverse Gaussian'' approximation is important in many applications, in
particular (see Section~\ref{sdyjumt}) in the problem of finding a solution of
non-linear equation \eqref{dsfgnhrgmn}; in the risk context this solution is
called ``non-ruin capital''; it is used to construct different controls in the
models of long-term insurance solvency regulation (see, e.g.,
\cite{Malinovskii-2015-b}--\cite{Malinovskii-2015-a}).

\end{document}